\documentclass{amsart}
\usepackage{amsmath}
\usepackage{amsthm}
\usepackage{amssymb}
\usepackage{graphicx}

\theoremstyle{plain}
\newtheorem{thm}{Theorem}[section]
\newtheorem{prop}[thm]{Proposition}

\newtheorem{lem}[thm]{Lemma}
\theoremstyle{remark}

\theoremstyle{definition}
\newtheorem{ex}[thm]{Example}
\newtheorem{defn}[thm]{Definition}
\numberwithin{equation}{section}
\begin{document}

\title[\lambdall]{Applications of a finite-dimensional duality\\
principle to some prediction problems}

\author[\lambdall]{Yukio Kasahara, Mohsen Pourahmadi and Akihiko Inoue}

\address{Department of Mathematics \\
Hokkaido University \\
Sapporo 060-0810, Japan}
\email{y-kasa@math.sci.hokudai.ac.jp}

\address{Division of Statistics\\
Northern Illinois University\\
DeKalb, IL 60115-2854, USA}
\email{pourahm@math.niu.edu}

\address{Department of Mathematics \\
Hokkaido University \\
Sapporo 060-0810, Japan}
\email{inoue@math.sci.hokudai.ac.jp}

\date{August 29, 2007}


\newcommand{\abs}[1]{\lvert#1\rvert}
\newcommand{\nor}[1]{\lVert#1\rVert}
\newcommand{\Na}{\mathbb{N}}
\newcommand{\Z}{\mathbb{Z}}
\newcommand{\R}{\mathbb{R}}
\newcommand{\C}{\mathbb{C}}
\newcommand{\M}{{\mathcal{M}}}
\newcommand{\Hc}{{\mathcal{H}}}
\newcommand{\Pc}{{\mathcal{P}}}
\newcommand{\ba}{\mbox{\boldmath$a$}}
\newcommand{\bob}{\mbox{\boldmath$b$}}
\newcommand{\bc}{\mbox{\boldmath$c$}}
\newcommand{\bd}{\mbox{\boldmath$d$}}
\newcommand{\be}{\mbox{\boldmath$e$}}
\newcommand{\bal}{\mbox{\boldmath$\alpha$}}
\newcommand{\bgam}{\mbox{\boldmath$\gamma$}}
\newcommand{\dsum}{\displaystyle\sum}

\newcounter{romnum}
\newenvironment{romlist}{\begin{list}{\roman{romnum}.}
{\usecounter{romnum}\setlength{\topsep}{1pt}
\setlength{\itemsep}{1pt}}
\rm}{\end{list}}


\begin{abstract}
Some of the most important results in prediction theory and time
series analysis when finitely many values are removed from or
added to its infinite past have been obtained using difficult and
diverse techniques ranging from duality in Hilbert spaces of
analytic functions (Nakazi, 1984) to linear regression in
statistics (Box and Tiao, 1975). 
We unify these results via a finite-dimensional duality lemma and
elementary ideas from the linear algebra. 
The approach reveals the inherent finite-dimensional character of
many difficult prediction problems, the role of duality and 
biorthogonality for a finite set of random variables. 
The lemma is particularly useful when the number of missing values
is small, like one or two, as in the case of Kolmogorov and Nakazi
prediction problems.
The stationarity of the underlying process is not a requirement. 
It opens up the possibility of extending such results to
nonstationary processes.
\end{abstract}

\maketitle

\renewcommand{\thefootnote}{\fnsymbol{footnote}} \footnote[0]{%
2000\textit{\ Mathematics Subject Classification} Primary 62M20;
Secondary 60G10; 60G25.}

\section{Introduction}

Irregular observations, missing values and outliers are common in
time series data
(Box and Tiao (1975), Brubacher and Wilson (1976)).
A framework for dealing with such anomalies is that of
$X=\{X_t\}_{t\in\Z}$ being a $\C$-valued, mean-zero, weakly
stationary stochastic process with the autocovariance function
$\gamma=\{\gamma_k\}_{k\in\Z}$ and the spectral density function
$f$: 
$E[X_k\bar{X_l}]=\gamma_{k-l}
=(2\pi)^{-1}\int_{-\pi}^{\pi}
e^{-i(k-l)\lambda}f(\lambda)d\lambda$. 
Then, the problem can be formulated as that of predicting or
approximating an unknown value $X_0$ based on the observed values
$\{X_t;t\in S\}$ for a given index set $S\subset\Z\setminus\{0\}$
and the knowledge of the autocovariance of the process.
Such a problem is quite important to applications in business,
economics, engineering, physical and natural sciences etc., and
belongs to the area of prediction theory of stationary stochastic
processes developed by Wiener (1949) and Kolmogorov (1941)
(see also Pourahmadi (2001)).
By restricting attention to linear predictors and using the 
least-squares criterion to assess the goodness of predictors, a
successful solution seeks to address the following two goals:
\begin{itemize}
\item[(P$_1$)]
Express the linear least-squares predictor of $X_0$, 
denoted by $\hat X_0(S)$, and the prediction error 
$X_0-\hat X_0(S)$ in terms of the observable $\{X_t;t\in S\}$.
\item[(P$_2$)]
Express the prediction error variance
$\sigma^2(S)=\sigma^2(f,S):=E\abs{X_0-\hat X_0(S)}^2$ in terms of
$f$.
\end{itemize}

The link between solutions of finite and infinite past prediction
problems serves as a natural bridge between time series analysis
and prediction theory.
From the dawn of modern time series analysis, the works of
Slutsky and Yule in the 1920's and Wold in the 1930's have been 
instrumental in achieving the goal (P$_1$) in the time-domain
using the {\it finite\/} past.
Subsequently, the classes of autoregressive (AR), moving-average
(MA) and mixed autoregressive and moving-average (ARMA) models
have played major roles in the development of time-domain
techniques using the autocovariance function of the process 
(see Box {\it et al.}\ (1994)). 
Nowadays, these techniques are implemented by solving the
Yule--Walker equations via the celebrated Durbin--Levinson
algorithm and the innovation algorithm
(see Brockwell and Davis (1991)). 
On the other hand, the spectral-domain techniques in prediction of
stationary processes, advocated by Kolmogorov and Wiener in the
early 1940's, rely on the spectral representations of the process
and its covariance
(Kolmogorov (1941), Wiener (1949), Pourahmadi (2001)).

The focus in prediction theory is more on the goal (P$_2$).
The celebrated Szeg\"o--Kolmogorov--Wiener theorem gives the
variance of the one-step ahead prediction error based on the
{\em infinite\/} past indexed by the ``half-line" 
$S_0:=\{ \dots,-2,-1\}$ by
\begin{equation}
\sigma^2(f,S_0)=\exp\left(\frac1{2\pi}\int_{-\pi}^\pi
\log f(\lambda)d\lambda\right)>0
\label{eq:skw}
\end{equation}
if $\log f$ is integrable, and otherwise $\sigma^2(S_0)=0$. 
However, when the first $n$ consecutive integers are removed from
$S_0$ or for the index set $S_{-n}:=\{\ldots,-n-2,-n-1\}$,
$n\ge 0$, the formula for the $(n+1)$-step prediction error
variance (Wold (1938), Kolmogorov (1941)) is 
\begin{equation}
\sigma^2(f,S_{-n})
=\vert b_0\vert^2+\vert b_1\vert^2+\cdots+\vert b_n\vert^2,
\qquad n=0,1,\ldots,
\label{eq:kw}
\end{equation}
where $\{b_j\}$, the MA coefficients of the process, is related to
the Fourier coefficients of $\log f$ and
$\abs{b_0}^2=\sigma^2(S_0)$
(see Nakazi and Takahashi (1980) and Pourahmadi (1984);
see also Section \ref{sec:FPprob} below).

A result similar to (\ref{eq:skw}) for the interpolation of a
single missing value corresponding to the index set
$S_\infty:=\Z\setminus\{0\}$ was obtained by Kolmogorov (1941).
Specifically, the interpolation error variance is given by
\begin{equation}
\sigma^2(f,S_\infty)
=\left(\frac1{2\pi}\int_{-\pi}^\pi
f(\lambda)^{-1}d\lambda\right)^{-1}>0
\label{eq:s}
\end{equation}
if $f^{-1}\in L^1:=L^1([-\pi,\pi], d\lambda/(2\pi))$, and
otherwise $\sigma^2(S_\infty)=0$. 
The corresponding prediction problem for the smaller index set 
$S_n:=\{\ldots,n-1,n\}\setminus\{0\}$, $n\ge 0$, was stated as
open in Rozanov (1967, p.~107) and is perhaps one of the most
challenging problems in prediction theory next to (\ref{eq:skw}). 
The index set $S_n$ is, indeed, of special interest as it forms a
bridge connecting $S_0$ and $S_\infty$;
it reduces to $S_0$ when $n=0$ and tends to $S_\infty$ as
$n\to\infty$.
In a remarkable paper in 1984, Nakazi using delicate, but
complicated analytical techniques (and assuming that
$f^{-1}\in L^1$) showed that
\begin{equation}
\sigma^2(f,S_n)=\left(\vert a_0\vert^2+\vert a_1\vert^2
+\cdots+\vert a_n\vert^2\right)^{-1},
\qquad n=0,1,\ldots,
\label{eq:n}
\end{equation}
where $\{a_j\}$ is related to the AR parameters of the process 
(see Section 3 below).

From (\ref{eq:kw}) and (\ref{eq:n}), the question naturally arises
as why there is such an ``inverse-dual'' relationship between
them. 
In this regard, it is worth noting that Nakazi's technique, if
interpreted properly, amounts to reducing computation of
$\sigma^2(f,S_n)$ to that of the $(n+1)$-step prediction error 
variance of another stationary process $\{Y_t\}$ with the
spectral density function $f^{-1}$ which turns out to be the
{\em dual\/} of $\{ X_t \}$ (see Definition 2.1 and Section 3.5). 
His result and technique have spawned considerable research in
this area in the last two decades; 
see Miamee and Pourahmadi (1988), Miamee (1993), 
Cheng {\em et al.}~(1998), Frank and Klotz (2002), 
Klotz and Riedel (2002) and Bondon (2002). 
A unifying feature of most of the known results thus far seems to
be a fundamental duality principle (Cheng {\it et al.}\ (1998),
Urbanik (2000)) of the form
\begin{equation}
\sigma^2(f,S)\cdot\sigma^2(f^{-1},S^c)=1,
\label{eq:dual}
\end{equation}
where $S^c$ is the complement of $S$ in $\Z\setminus\{0\}$ and 
$f^{-1}\in L^1$. 
The first occurrence of (\ref{eq:dual}) seems to be in the 1949
Russian version of Yaglom (1963) for the case of deleting finitely
many points from $S_\infty$. 
Proof of (\ref{eq:dual}), in general, like those of the main
results in Nakazi (1984), Miamee and Pourahmadi (1988),
Cheng {\it et al.}\ (1998), and Urbanik (2000), is long,
unintuitive and relies on duality techniques from functional and
harmonic analysis and requires $f^{-1}\in L^1$ which is not
natural for an index set like $S_n$.
Surprisingly, a version of (\ref{eq:dual}) in a rather disguised
form was developed in Grenander and Rosenblatt (1954, Theorem 1), 
as the limit of a quadratic form involving Szeg\"o's orthogonal
polynomials on the unit circle, see also Simon (2005, p.165). 
Unfortunately, it had remained dormant and not used in the context
of prediction theory, except in Pourahmadi (1993).

In this paper, we establish a finite-dimensional duality principle
(Lemma \ref{L:int}), which encapsulates (\ref{eq:dual}) in a
transparent and useful manner. 
The concept of dual of a random vector plays a central role as
does the Cholesky decomposition of its covariance matrix. 
We use this duality principle to unify and solve some prediction 
problems related to removing a finite number of indices from 
$S_n$ and $S_{\infty}$. 
The outline of the paper is as follows. 
In Section \ref{sec:MDL}, we present the main lemma, some
auxiliary facts about dual of a random vector and their
consequences for computing the prediction error variances and
predictors. 
In Section \ref{sec:FPprob}, using the lemma we first solve three
finite prediction problems for $X_0$ based on the knowledge of
$\{X_t; t\in K\}$ with $K = \{-m,\dots,n\}\setminus (M\cup\{0\})$,
$m, n \geq 0$, where $M$, the index set of the missing values, is
relatively small.
Then we obtain the solutions of Kolmogorov, Nakazi, and Yaglom's
prediction problems in a unified manner by studying the limit of
the solutions by letting $m\to\infty$, followed by $n\to\infty$.
In particular, we find explicit formula for the dual of the
process $\{X_t; t\le n\}$ for a fixed $n$, which does not seem to
be possible using the technique of Urbanik (2000), Klotz and
Riedel (2002) and Frank and Klotz (2002).
This is useful in developing series representations for predictors
and interpolators, and sheds light on the approaches of
Bondon (2002) and Salehi (1979).
In Section \ref{sec:discu}, we close the paper with some
discussions.

Finally, we should point out that the two simple formulas
(\ref{eq:kw}) and (\ref{eq:n}) and their extensions provide
explicit and informative expressions for the prediction error
variances.
Like their predecessors (\ref{eq:skw}) and (\ref{eq:s}), they
serve as yardsticks to assess the impact (worth) of observations
in predicting $X_0$ when they are added to or deleted from the
infinite past and highlight the role of the autoregressive and
moving-average parameters for this purpose; 
see Pourahmadi and Soofi (2000).
In fact, Bondon (2002, Theorem 3.3; 2005) shows that a finite
number of missing values do not affect the prediction of $X_0$ if
and only if the AR parameters corresponding to the indices of
those missing values are zero.
Furthermore, the examples in Section \ref{sec:FPprob} indicate how
the interpolators of the missing values can be computed rigorously
without resorting to formal derivations (Box and Tiao (1975),
Brubacher and Wilson (1976) and Budinsky (1989)).


\section{A Finite-Dimensional Duality Principle} \label{sec:MDL}

In this section, an elementary result is stated as a
finite-dimensional duality lemma, which we use in Section 3 to
solve and unify various challenging prediction problems through
the limit of the solutions of their finite past counterparts.

For a finite index set $N$, let $H_N$ be the class of vectors
$X=(X_j)_{j\in N}$ of random variables with zero-mean and finite
variance on a probability space $(\Omega,\mathcal F,P)$:
\[
H_N:=\{X=(X_j)_{j\in N};
\ X_j\in L^2(\Omega,\mathcal{F},P),\ E[X_j]=0,\ j\in N\}.
\]
As usual, we consider the inner product $(Y,Z):=E[Y\bar{Z}]$ and
norm $\nor{Y}:=E[\vert Y\vert^2]^{1/2}$ for random variables
in $L^2(\Omega,\mathcal F,P)$. 

\begin{defn}
Let $N$ be a finite index set and $X\in H_N$. 
A random vector $Y\in H_N$ is called the {\it dual of $X$} if it
satisfies the following conditions:
\begin{itemize}
\item[{\rm (i)}]
The components $Y_j$, $j\in N$, belong to
$\mathrm{sp}\{X_k;k\in N\}$.
\item[{\rm (ii)}] $X$ and $Y$ are {\em biorthogonal\/}:
$(X_i,Y_j)=\delta_{ij}$ for  $i,j\in N,$ or
$\mathrm{Cov}(X,Y)=I$. 
\end{itemize}
\end{defn}

For $X\in H_N$, $l\in N$ and $K\subset N$, we write $\hat X_l(K)$
for the linear least squares predictor of $X_l$ based on
$\{X_k; k\in K\}$, i.e., the orthogonal projection of $X_l$ onto
$\mathrm{sp}\{X_k; k\in K\}$. 
For the sake of completeness and ease of reference, in the next
two propositions we summarize the characterization, interpretation
and other basic information about the dual of a random vector in
terms of its covariance matrix and certain prediction errors.

\begin{prop}\label{prop:3.1aaa}
Let $N$ be a finite index set and $X\in H_N$. 
Then, the following conditions are equivalent:
\begin{itemize}
\item[{\rm (1)}]
The components $X_j$, $j\in N$, of $X$ are linearly independent.
\item[{\rm (2)}]
The covariance matrix $\Gamma=(\gamma_{i,j})_{i,j\in N}$ 
of $X$ with $\gamma_{i,j}=(X_i,X_j)$ is nonsingular.
\item[{\rm (3)}]
$X$ is minimal:
$X_j\notin \mathrm{sp} \{X_i; i\in N,\ i\neq j\}$ for $j\in N$.
\item[{\rm (4)}]
$X$ has a dual.
\end{itemize}
\end{prop}

\begin{proof}
Clearly, (1)--(3) are equivalent. 
Assume (3) and define $Y=(Y_j)_{j\in N}\in H_N$ by 
$Y_j=(X_j - \hat X_j(N_j))/
\nor{X_j-\hat X_j(N_j)}^2$, where $N_j:=N\setminus\{j\}$. 
Then $Y_j$ belongs to $\mathrm{sp} \{X_k; k\in N\}$,
and $(X_i,Y_j)=\delta_{ij}$ holds:
\[
(X_j,Y_j)
=\frac{(X_j,X_j-\hat X_j(N_j))}
{\nor{X_j-\hat X_j(N_j)}^2}
=\frac{(X_j-\hat X_j(N_j),X_j-\hat X_j(N_j))}
{\nor{X_j-\hat X_j(N_j)}^2}
=1,
\]
and for $i\neq j$,
\[
(X_i,Y_j)
=\frac{(X_i,X_j-\hat X_j(N_j))}
{\nor{X_j-\hat X_j(N_j)}^2}
=0.
\]
Thus $Y$ is a dual of $X$, and hence (4). 
Conversely, assume (4) and let $Y$ be a dual of $X$. 
If $X$ is not minimal, then there exists $j\in N$ such that 
$X_j\in\mathrm{sp}\{X_i; i\in N,\ i\neq j\}$, that is,
$X_j=\sum_{i\ne j}c_iX_i$ for some $c_i\in\C$, and, since
$(X_i,Y_j)=0$ for $i\ne j$, we have
$(X_j,Y_j)=\sum_{i\ne j}c_i(X_i,Y_j)=0$.
However, this contradicts $(X_j,Y_j)=1$. 
Thus, $X$ is minimal, and (3) follows. 
\end{proof}

The proof reveals the importance of the ``standardized"
interpolation errors of components of $X$ in defining its dual.
More explicit representations and other properties of the dual
are given next.

\begin{prop}\label{prop:3.2aaa}
For a finite index set $N$, let $X\in H_N$ 
with covariance matrix $\Gamma$. 
Assume that $X$ has a dual $Y$.
Then the following assertions hold:
\begin{itemize}
\item[{\rm (1)}] The dual $Y$ is unique.
\item[{\rm (2)}] The dual $Y$ is given by 
$Y_j=(X_j-\hat X_j(N_j))
/\Vert X_j-\hat X_j(N_j)\Vert^2$ with 
$N_j:=N\setminus\{j\}$ for $j\in N$. 
\item[{\rm (3)}] 
The dual $Y$ is also given by $Y=\Gamma^{-1}X$ 
or $Y_i=\sum_{j\in N}\gamma^{i,j}X_j$, $i\in N$, where  
$\Gamma^{-1}=(\gamma^{i,j})_{i,j\in N}$. 
\item[{\rm (4)}]The covariance matrix of $Y$ is equal to
$\Gamma^{-1}$.
\item[{\rm (5)}] The dual of $Y$ is $X$.
\item[{\rm (6)}] 
$\mathrm{sp} \{X_j; j\in N\}=\mathrm{sp} \{Y_j; j\in N\}$.
\end{itemize}
\end{prop}

\begin{proof}
First, we prove (1). 
Let $Z$ be another dual of $X$ and $j\in N$ be fixed. 
Then $(X_i,Y_j-Z_j)=0$ for all $i\in N$.
However, since $Y_j-Z_j\in\mathrm{sp}\{X_k; k\in N\}$, it follows
that $Y_j=Z_j$ and hence (1). 
(2) follows from the proof of Proposition \ref{prop:3.1aaa}. 
To prove (3) and (4), we put $Y=\Gamma^{-1}X$. 
Then $Y_j\in\mathrm{sp}\{X_k; k\in N\}$.
Since $\Gamma^{-1}$ is Hermitian, we have
\[
\begin{split}
&\mathrm{Cov}(X,Y)
=\mathrm{Cov}(X,X)\,\Gamma^{-1}=\Gamma\Gamma^{-1}=I,\\
&\mathrm{Cov}(Y,Y)
=\Gamma^{-1}\mathrm{Cov}(X,X)\,\Gamma^{-1}
=\Gamma^{-1}\Gamma\,\Gamma^{-1}
=\Gamma^{-1}.
\end{split}
\]
Thus (3) and (4) follow. 
Finally, we obtain (5) and (6) from (3) and (4). 
\end{proof}

From the two representations in Proposition \ref{prop:3.2aaa} (2),
(3) for the dual $Y$, we find the following representation for the
standardized interpolation error:
\[
\frac{X_i-\hat X_i(N_i)}
{\nor{X_i-\hat X_i(N_i)}^2}
=\sum\nolimits_{j\in N}\gamma^{i,j}X_j\quad
\mbox{with}\quad N_i=N\setminus\{i\}. 
\]
In particular, $\gamma^{i,i}=1/\nor{X_i-\hat X_i(N_i)}^2$. 
Notice that these equalities hold even if $\Gamma$ is not a
Toeplitz matrix or $X$ is not a segment of a stationary process. 
For some statistical/physical interpretations of the entries of
$\Gamma^{-1}$, the inverse of a stationary covariance matrix, see
Bhansali (1990) and references therein.

Now, we are ready to state the main duality lemma.

\begin{lem}\label{L:int}
Let $N$ be a finite index set. 
Assume that $X\in H_N$ has the dual $Y\in H_N$ 
and that 
$K$, $M$ and a singleton $\{l\}$ partition $N$, i.e.,
\[
N=K\cup\{l\}\cup M\qquad\mbox{\rm (disjoint union)}.
\]
Then the following equalities hold:
\begin{itemize}
\item[{\rm (a)}]
$\displaystyle 
X_l-\hat X_l(K)
=\frac{Y_l-\hat Y_l(M)}{\Vert Y_l-\hat Y_l(M)\Vert^2}$.
\item[{\em (b)}]
$\displaystyle\Vert X_l - \hat X_l(K)\Vert
=\frac{1}{\Vert Y_l-\hat Y_l(M)\Vert}$.
\end{itemize}
\end{lem}

\begin{proof}
Since $X$ and $Y$ are minimal and biorthogonal, $X_l-\hat X_l(K)$
and $Y_l-\hat Y_l(M)$ are nonzero and belong to the same
one-dimensional space, that is, the orthogonal complement of
$\mathrm{sp}\{X_j;j\in K\}\oplus\mathrm{sp}\{Y_j;j\in M\}$ in
$\mathrm{sp}\{X_j;j\in N\}$.
Therefore, one is a multiple of the other; for some $c\in\C$, 
\[
X_l-\hat X_l(K)
=c\,
\frac{Y_l-\hat Y_l(M)}{\nor{Y_l-\hat Y_l(M)}^2}.
\]
But, since $c$ is equal to
\[
c\,\frac{(Y_l-\hat Y_l(M),Y_l-\hat Y_l(M))}
{\nor{Y_l-\hat Y_l(M)}^2}
=c\,\frac{(Y_l-\hat Y_l(M),Y_l)}
{\nor{Y_l-\hat Y_l(M)}^2}
=(X_l-\hat X_l(K),Y_l)
=(X_l,Y_l)=1,
\]
we get (a) and (b) and hence the lemma. 
\end{proof}

In the applications in Section 3, we use this duality in the form
of the next lemma which gives a way of computing the predictor
coefficients and prediction error variance using the inverse
matrix $\Gamma^{-1}=(\gamma^{i,j})$.

\begin{lem}\label{L:snd}
Let $N$, $X=(X_j)_{j\in N}$, $Y=(Y_j)_{j\in N}$, $K$, $M$ and
$\{l\}$ be as in Lemma \ref{L:int} with
$\Gamma=(\gamma_{i,j})_{i,j\in N}$ the covariance matrix of $X$
and $\Gamma^{-1}=(\gamma^{i,j})_{i,j\in N}$. 
Then
\begin{align}
&X_l-\hat X_l(K)
=\sum\nolimits_{i\in M\cup\{l\}}\alpha'_iY_i,
\label{eq:fff}\\
&\Vert X_l-\hat X_l(K)\Vert^2=\alpha'_l,
\label{eq:ggg}
\end{align}
where $(\alpha'_i)_{i\in M\cup\{l\}}$ is the solution to 
the following system of linear equations:
\begin{equation}
\sum\nolimits_{i\in M\cup\{l\}}\alpha'_i\gamma^{i,j}
=\delta_{lj}, 
\qquad j\in M\cup\{l\}.
\label{eq:hhh}
\end{equation}
In particular, 
the prediction error
variance $\sigma_l^2(K)=\Vert X_l-\hat X_l(K)\Vert^2$ is given by
\begin{equation}
\sigma_l^2(K)
=\mbox{the $(l,l)$-entry of the inverse of 
$(\gamma^{i,j})_{i,j\in M\cup\{l\}}$},
\label{eq:jjj}
\end{equation}
and the predictor coefficients $\alpha_k$ in 
$\hat X_l(K)=\sum_{k\in K}\alpha_kX_k$ are given by 
\begin{equation}
\alpha_k=-\sum\nolimits_{i\in M\cup\{l\}}\alpha'_i\gamma^{i,k},
\qquad k\in K,
\label{eq:iii}
\end{equation}
whence we have
\begin{equation}
\hat X_l(K)=-\sum\nolimits_{k\in K}
\left(\sum\nolimits_{i\in M\cup\{l\}}\alpha'_i\gamma^{i,k}\right)
X_k.
\label{eq:iv}
\end{equation}
\end{lem}

\begin{proof}
Since $Y_j$'s are linearly independent, Lemma \ref{L:int} (a)
shows that $X_l-\hat X_l(K)$ is uniquely expressed in the form
(\ref{eq:fff}).
Then $\alpha'_l=\Vert Y_l-\hat Y_l(M)\Vert^{-2}$, which,
in view of Lemma \ref{L:int} (b), is equal to
$\Vert X_l-\hat X_l(K)\Vert^2$, and (\ref{eq:ggg}) holds.
Since $(X_i,Y_j)=\delta_{ij}$ and $(Y_i,Y_j)=\gamma^{i,j}$, the
predictor coefficients
$\alpha_k$ in $\hat X_l(K)=\sum_{k\in K}\alpha_kX_k$ satisfy
\[
\alpha_k
=(\hat X_l(K),Y_k)
=\left(X_l-\sum\nolimits_{i\in M\cup\{0\}}\alpha'_iY_i,
Y_k\right)
=-\sum\nolimits_{i\in M\cup\{0\}}\alpha'_i\gamma^{i,k}.
\]
Thus (\ref{eq:iii}), whence (\ref{eq:iv}).
Similarly, for $j\in M\cup\{l\}$, we have $(\hat X_l(K),Y_j)=0$
and
\[
\sum\nolimits_{i\in M\cup\{l\}}\alpha'_i\gamma^{i,j}
=\sum\nolimits_{i\in M\cup\{l\}}\alpha_i'(Y_i,Y_j)
=(X_l-\hat X_l(K),Y_j)=\delta_{lj}.
\]
Therefore, (\ref{eq:hhh}) follows.
Finally, we obtain (\ref{eq:jjj}) from (\ref{eq:ggg}) and
(\ref{eq:hhh}).
\end{proof}

Recall that the predictor coefficients $\alpha_k=\alpha_{k,l}(K)$
in $\hat X_l(K)=\sum_{k\in K}\alpha_kX_k$ and the prediction error
variance $\sigma^2=\sigma_l^2(K)=\Vert X_l-\hat X_l(K)\Vert^2$ are
traditionally computed from $(\gamma_{i,j})_{i,j\in K\cup\{l\}}$
by solving the normal equations:
\begin{equation}
\left\{
\begin{aligned}
&\sum\nolimits_{k\in K}\alpha_{k}\gamma_{k,j}=\gamma_{l,j},
\qquad j\in K,\\
&\sigma^2
=\gamma_{l,l}-\sum\nolimits_{k\in K}\alpha_{k}\gamma_{k,l}.
\end{aligned}
\right.
\label{eq:normal}
\end{equation}
Alternatively, one could write the above as an analogue of the 
Yule--Walker equations:
\begin{equation}
\gamma_{l,j} - \sum\nolimits_{k\in K}\alpha_{k}\gamma_{k,j}
=\delta_{lj}\sigma^2,
\qquad 
j\in K\cup\{l\}.
\label{eq:YW-eq}
\end{equation}
Then $\sigma^2=\sigma_l^2(K)$ can be identified as
\begin{equation}
\sigma_l^2(K)
=\left[\,\mbox{the $(l,l)$-entry of the inverse of 
$(\gamma_{i,j})_{i,j\in K\cup\{l\}}$}\right]^{-1}.
\label{eq:sigma-minus2}
\end{equation}
In addition, using the Cramer's rule, one may write $\sigma^2$ in
(\ref{eq:sigma-minus2}) as the ratio of the two relevant
determinants: 
$\sigma^2=\det(\gamma_{i,j})_{i,j\in K\cup\{l\}}/
\det(\gamma_{i,j})_{i,j\in K}$.

In spite of the simplicity of
(\ref{eq:normal})--(\ref{eq:sigma-minus2}), they are not
convenient for the study of the asymptotic behaviors of the
predictor coefficients and predictor variance as $K$ gets large.
The method of computation in Lemma \ref{L:snd} becomes
particularly useful when $K$ is large but $M$ is small
(see Section 3.2 below).


\section{Applications to Prediction Problems} \label{sec:FPprob}

In this section, we illustrate the role of the finite duality 
principle (Lemmas \ref{L:int} and \ref{L:snd}) 
in unifying some diverse prediction problems
for a zero-mean, weakly stationary process $\{X_j\}_{j\in\Z}$ 
with the autocovariance function $\gamma=\{\gamma_j\}_{j\in\Z}$:
$\gamma_{i-j}=(X_i,X_j)$.

For simplicity, we assume that $\{X_j\}_{j\in\Z}$ is purely
nondeterministic, so it admits the MA representation
(Wold decomposition)
\begin{equation}
X_j=\sum\nolimits_{k=-\infty}^j b_{j-k}\varepsilon_k, 
\qquad
j\in\Z,
\label{eq:MA}
\end{equation}
where $\{\varepsilon_j\}_{j\in\Z}$ is the normalized innovation of
$\{X_j\}_{j\in\Z}$ defined by
\[
\varepsilon_j
:=\{X_j-\hat X_j(\{\dots,j-2,,j-1\})\}/
\Vert X_j-\hat X_j(\{\dots,j-2,,j-1\})\Vert,
\qquad j\in\Z,
\]
and $\{b_k\}_{k=0}^\infty$ is the MA coefficients 
given by $b_k:=(X_0,\varepsilon_{-k})$. 
We define a sequence of complex numbers $\{a_k\}_{k=0}^\infty$ by
the relation
\begin{equation}
\sum\nolimits_{k=0}^jb_ka_{j-k}=\delta_{0j},
\qquad
j\geq 0.
\label{eq:ba}
\end{equation}
If the series $\sum_{j=0}^\infty a_jX_{-j}$ is mean-convergent,
then (\ref{eq:MA}) can inverted as
\begin{equation}
\varepsilon_j
=\sum\nolimits _{k=-\infty}^ja_{j-k}X_k,
\qquad
j\in\Z.
\label{eq:AR}
\end{equation}
This is essentially the same as the AR representation (see
Pourahmadi (2001)), and we call $\{a_k\}$ the AR coefficients of
$\{X_j\}_{j\in\Z}$.
As suggested in (\ref{eq:kw}) and (\ref{eq:n}), these $\{b_k\}$
and $\{a_k\}$ play an important role in prediction problems.


\subsection{Finite Prediction Problems with Missing Values}
\label{sec:4.2}

Let $M$ be a finite set of integers that does not contain zero.
Throughout this section, it represents the index set of missing
(unknown) values when predicting $X_0$.
For given $M$, we take the integers $m,n\ge 0$ so large that 
$M\subset N:=\{-m,\dots,n\}$, and put
$K=N\setminus(M\cup\{0\})$, which represents the index set of the
observed values, so that we have
the partition $N=K\cup\{0\}\cup M$ as in Lemma \ref{L:int}.
We start with the prediction problem for a finite index set $K$. 
Once the problem is solved for such a $K$, the solutions for
infinite index sets $S_n\setminus M$ and $S_{\infty}\setminus M$ 
are obtained by taking the limit of the solutions, first as
$m\to\infty$, and then $n\to\infty$.

Traditionally, the coefficients of the finite linear predictor
$\hat X_0(K)$ and its prediction error variance
$\sigma^2(K)=\nor{X_0-\hat X_0(K)}^2$ are expressed in terms of
the covariance function $\gamma$, using the normal equations
(\ref{eq:normal}). 
However, the results so obtained are not convenient for studying
the asymptotic behaviors of the predictor coefficients as
$m\to\infty$ and/or $n\to\infty$. 
The problem can be made much simpler by the finite duality 
principle and some fundamental facts about the finite MA and AR
representations, as we explain now (see also Pourahmadi (2001)).

For the future segment $\{X_j\}_{j=0}^\infty$ of the process, we
define its normalized innovation
$\{\varepsilon_{j,0}\}_{j=0}^\infty$ by the Gram--Schmidt method:
$\varepsilon_{0,0}:=X_0/\Vert X_0\Vert$ and 
\[
\varepsilon_{j,0}
:=\{X_j-\hat X_j(\{0,\dots,j-1\})/
\Vert X_j-\hat X_j(\{0,\dots,j-1\})\Vert,
\qquad j\ge 1.
\]
Then $\{X_j\}$ and $\{\varepsilon_{j,0}\}$ admit the following
finite MA and AR representations:
\[
X_j=\sum\nolimits_{k=0}^jb_{j-k,j}\varepsilon_{k,0},
\qquad
\varepsilon_{j,0}=\sum\nolimits_{k=0}^ja_{j-k,j}X_k,
\qquad
j\geq0.
\]
Here $\{b_{k,j}\}_{k=0}^j$ is defined by
$b_{k,j}:=(X_j,\varepsilon_{j-k,0})$ and $\{a_{k,j}\}_{k=0}^j$ by
\[
\sum\nolimits_{k=i}^jb_{j-k,j}a_{k-i,k}=\delta_{ij}
\qquad\mbox{or}\qquad
\sum\nolimits_{k=i}^ja_{j-k,j}b_{k-i,k}=\delta_{ij},
\qquad
i\leq j.
\]
These finite MA and AR coefficients converge to their infinite
counterparts:
\begin{equation}
\lim_{j\to\infty}b_{k,j}=b_k,
\qquad
\lim_{j\to\infty}a_{k,j}=a_k.
\label{eq:integ}
\end{equation}
If we consider $\{X_j\}_{j=-m}^\infty$ instead of
$\{X_j\}_{j=0}^\infty$, then by stationarity, it follows that
\begin{equation}
X_j
=\sum\nolimits_{k=-m}^jb_{j-k,{m+j}}\varepsilon_{k,-m},
\quad
\varepsilon_{j,-m}
=\sum\nolimits_{k=-m}^ja_{j-k,{m+j}}X_k,
\quad
j\geq-m,
\label{eq:4.2aaa}
\end{equation}
where $\{\varepsilon_{j,-m}\}_{j=-m}^\infty$ is the normalized
innovation of $\{X_j\}_{j=-m}^\infty$ defined in the same way.
We notice that
\begin{equation}
\varepsilon_j
=\lim_{m\to\infty}\varepsilon_{j,-m},
\qquad j\in\Z.
\label{eq:e-e}
\end{equation}
Thus, the representations in (\ref{eq:4.2aaa}) reduce to
(\ref{eq:MA}) and (\ref{eq:AR}) as $m\to\infty$.

Recall that $N=\{-m,\ldots,n\}$ and let $X$ be the vector
$(X_j)_{j\in N}$ with covariance matrix
$\Gamma=(\gamma_{i-j})_{i,j\in N}$. 
From Proposition \ref{prop:3.2aaa} (3), its dual $Y$ is given by
$Y=\Gamma^{-1}X$. 
Let $\varepsilon$ be the normalized innovation vector of $X$,
i.e., $\varepsilon:=(\varepsilon_{j,-m})_{j\in N}$. 
Then it follows from (\ref{eq:4.2aaa}) that
\[
X=B\varepsilon,\qquad \varepsilon=AX,
\]
where $A$ and $B$ are the lower triangular matrices with
$(i,j)$-entries $a_{i-j,m+i}$ and $b_{i-j,m+i}$ for
$-m\le j\le i\le n$, respectively.
Since $A=B^{-1}$ and $\Gamma=BB^*$, we have  
\[
\Gamma^{-1}=A^*A,
\qquad
Y=A^\ast\varepsilon.
\]
Thus, the $(i,j)$-entry $\gamma^{i,j}$ of $\Gamma^{-1}$ and the
$j$-th entry $Y_j$ of $Y$ have the representations
\begin{equation}
\gamma^{i,j}
=\sum\nolimits_{k=i\vee j}^n \bar{a}_{k-i,m+k}a_{k-j,m+k},
\qquad 
Y_j
=\sum\nolimits_{k=j}^{n}\bar{a}_{k-j,m+k}\varepsilon_{k,-m},
\label{eq:4.4aaa}
\end{equation}
which are certainly more conducive to studying their limits as
first $m\to\infty$ and then $n\to\infty$, see (\ref{eq:integ}) and
(\ref{eq:e-e}).

Now, we are ready to express the predictor $\hat X_0(K)$, the
prediction error $X_0-\hat X_0(K)$ and its variance $\sigma^2(K)$
as prescribed
by Lemma \ref{L:snd}.
In particular, it follows from (\ref{eq:fff}), (\ref{eq:jjj}) and
(\ref{eq:4.4aaa}) that 
\begin{equation}
\sigma^2(K)
=\mbox{the $(0,0)$-entry of the inverse of }
\left(\sum\nolimits_{k=i\vee j}^n\bar a_{k-i,m+k}a_{k-j,m+k}
\right)_{i,j\in M\cup\{0\}}
\label{eq:4.6aaa}
\end{equation}
and
\begin{equation}
X_0-\hat X_0(K)
=\sum\nolimits_{i\in M\cup\{0\}}\alpha_i^\prime
\left(\sum\nolimits_{k=i}^{n}\bar{a}_{k-i,m+k}\varepsilon_{k,-m}
\right),
\label{eq:4.5aaa}
\end{equation}
where $\alpha_i^\prime$'s are as in Lemma \ref{L:snd} with $l=0$.

To highlight some far-reaching consequences of (\ref{eq:4.6aaa})
and (\ref{eq:4.5aaa}), a few special cases corresponding to the
classical prediction problems of Kolmogorov (1941), Yaglom (1963)
and Nakazi (1984) are singled out and listed as examples in the
next section according to the cardinality of the index set $M$ of
the missing values.


\subsection{Examples} \label{sec:examp}

In this section, we discuss three distinct examples of the use of 
the finite duality principle and illustrate the process of
obtaining results for the two infinite index sets
$S=S_n\setminus M$ and $S=S_\infty\setminus M$.

Since $\{X_j\}_{j\in\Z}$ is purely nondeterministic, it has the
spectral density function $f$ with $\log f\in L^1$:
$\gamma_j
=(2\pi)^{-1}\int_{-\pi}^\pi e^{-ij\lambda}f(\lambda)d\lambda$.
Also, there exists an outer function $h$ in the Hardy class $H^2$
such that $f=\vert h\vert^2$ and $h(0)>0$, and we have
\begin{equation} \label{eq:follows}
h(z)=\sum\nolimits_{k=0}^\infty b_kz^k,
\qquad
\frac1{h(z)}=\sum\nolimits_{k=0}^\infty a_kz^k
\end{equation}
in the unit disc.
This shows that $f^{-1}\in L^1$ if and only if $\{a_k\}$ is square
summable.
Using (\ref{eq:follows}), which should be compared with
(\ref{eq:MA})--(\ref{eq:AR}), we can define the MA and AR
coefficients in an analytical way.

\begin{ex}[The Finite Kolmogorov--Nakazi Problem]
This is a finite interpolation problem corresponding to
$K=\{-m,\ldots,n\}\setminus\{0\}$ and $M=\phi$ (empty set),
and the solution of (\ref{eq:hhh}) is $\alpha'_0=1/\gamma^{0,0}$.
Consequently, from (\ref{eq:4.4aaa})--(\ref{eq:4.5aaa}), we have
\begin{equation}
\sigma^2(K)
=\left(\sum\nolimits_{k=0}^{n}\vert a_{k,m+k}\vert^2\right)^{-1}
\label{eq:4.8aaa}
\end{equation}
and 
\begin{equation}
X_0-\hat X_0(K)
=\left(\sum\nolimits_{k=0}^{n}\vert a_{k,m+k}\vert^2\right)^{-1}
\sum\nolimits_{k=0}^{n}\bar{a}_{k,m+k}\varepsilon_{k,-m}.
\label{eq:4.7aaa}
\end{equation}

Next, we show that (\ref{eq:4.8aaa}) and (\ref{eq:4.7aaa}) are
precursors of important results in prediction theory due to
Kolmogorov (1941), Masani (1960), and Nakazi (1984). 

The result (\ref{eq:n}) of Nakazi (1984) for
$S_n=\{\ldots,n-1,n\}\setminus\{0\}$ is obtained by taking the
limit of (\ref{eq:4.8aaa}) as $m\to\infty$
(without assuming $f^{-1}\in L^1$).
Indeed, by (\ref{eq:integ}), we see that (\ref{eq:4.8aaa}) gives
\begin{equation}
\sigma^2(S_n)
=\left(\sum\nolimits_{k=0}^{n}\vert a_{k}\vert^2\right)^{-1}.
\label{eq:4.11aaa}
\end{equation}
Also, in view of (\ref{eq:e-e}),
it follows from (\ref{eq:4.7aaa}) that
\begin{equation}
X_0-\hat X_0(S_n)
=\left(\sum\nolimits_{k=0}^{n}\vert a_{k}\vert^2\right)^{-1}
\sum\nolimits_{k=0}^{n}\bar{a}_{k}\varepsilon_{k}.
\label{eq:4.10aaa}
\end{equation}

The solution (\ref{eq:s}) of the Kolmogorov (1941) interpolation
problem with $S_\infty=\Z\setminus\{0\}$ follows from
(\ref{eq:4.11aaa}) by taking the limit as $n\to\infty$, provided
that $\{a_k \}$ is square summable.  
Thus, as in Kolmogorov (1941), assuming that $\{X_t\}$ is minimal
or $f^{-1}\in L^1$, we obtain
\[
\sigma^2\left(S_\infty \right)
=\left(\sum\nolimits_{k=0}^{\infty}\vert a_{k}\vert^2\right)^{-1}
=\left(\frac1{2\pi}
\int_{-\pi}^{\pi}f(\lambda)^{-1}d\lambda\right)^{-1}.
\]
Under the same minimality condition, the limit of
(\ref{eq:4.10aaa}) as $n\to\infty$, leads to 
\[
X_0-\hat X_0(S_\infty)
= \left(\sum\nolimits_{k=0}^{\infty}\vert a_{k}\vert^2\right)^{-1}
\sum\nolimits_{k=0}^{\infty}\bar{a}_{k}\varepsilon_{k},
\]
which is Masani's (1960) representation of the two-sided
innovation of $\{X_j\}$ at time $0$.
It is instructive to note that this is a moving average in terms
of the future innovations.
In fact, the source of such moving average representation can be
traced to (\ref{eq:4.4aaa}) and (\ref{eq:4.10aaa}).
A version of (\ref{eq:4.10aaa}) seems to have appeared first in
Box and Tiao (1975) in the context of intervention analysis;
see Pourahmadi (1989), and Pourahmadi (2001, Section 8.4) for a
more rigorous derivation, detailed discussion and connection with
outlier detection.
\end{ex}

Our second example corresponds to $M$ having cardinality one and
hence involves inversion of $2\times 2$ matrices, no matter how 
large $K$ is.

\begin{ex}[The Finite Past with a Single Missing Value]
This problem corresponds to $m>0$, $n=0$,
$K=\{-m,\dots,-1\}\setminus\{-u\}$ and  $M=\{-u\}$, where
$1\leq u\leq m$, so that $X_{-u}$ from the finite past of length
$m$ is missing. 
By (\ref{eq:4.4aaa}), the $2\times 2$ matrix for solving
(\ref{eq:hhh}) is
\[
\left( \begin{array}{cc}
\gamma^{-u,-u} & \gamma^{-u,0}\\
\gamma^{0,-u} & \gamma^{0,0}
\end{array} 
\right)
=\left( \begin{array}{cc}
\sum_{k=0}^u\vert a_{u-k, m- k}\vert^2 & a_{0,m}\,\bar{a}_{u,m}\\
\bar{a}_{0,m}\,a_{u,m} & \vert a_{0,m}\vert^2
\end{array} 
\right).
\]
Hence, using the subscript $m$ to emphasize the dependence on $m$,
we have
\begin{equation*}
\alpha'_{0,m}=\frac{1}{\Delta_m}
\sum\nolimits_{k=0}^u \vert a_{u-k,m-k}\vert^2,
\qquad 
\alpha'_{-u,m}=-\frac{\bar{a}_{0,m}a_{u,m}}{\Delta_m},
\end{equation*}
with the determinant
$\Delta_m
=\vert a_{0,m}\vert^2\sum_{k=1}^u \vert a_{u-k,m-k}\vert^2$. 
Thus, by (\ref{eq:4.6aaa}) and (\ref{eq:4.5aaa}),
\begin{equation}
\left\{ 
\begin{aligned}
&\sigma^2(K)
=\dfrac{\sum_{k=0}^u
\vert a_{u-k,m-k}\vert^2}{\vert a_{0,m}\vert^2
\sum_{k=1}^u \vert a_{u-k,m-k}\vert^2}, \\
&X_0-\hat X_0(K)
=\alpha'_{0,m}\bar{a}_{0,m}\varepsilon_{0,-m} 
+\alpha'_{-u,m}
\sum\nolimits_{k=0}^u\bar{a}_{u-k,m-k}\varepsilon_{-k,-m},
\end{aligned}
\right. 
\label{eq:ent}
\end{equation} 
and, taking the limit as $m\to\infty$, 
\begin{equation}
\left\{
\begin{aligned}
&\sigma^2(S_0\setminus\{-u\})
=\abs{b_0}^2
\dfrac{\sum_{k=0}^{u}\vert a_k\vert^2}
{\sum_{k=0}^{u-1}\vert a_k\vert^2},\\
&X_0-\hat X_0(S_0\setminus\{-u\})
=\alpha'_{0}\bar{a}_{0}\varepsilon_{0} 
+ \alpha'_{-u}\sum\nolimits_{k=0}^u\bar{a}_{u-k}\varepsilon_{-k},
\end{aligned}
\right.
\label{eq:4.15taking}
\end{equation}
where $\alpha'_0$ and $\alpha'_{-u}$ are the limits of
$\alpha'_{0,m}$ and $\alpha'_{-u,m}$, as $m\to\infty$,
respectively. 

The expressions in (\ref{eq:4.15taking}) were obtained first in
Pourahmadi (1992); see also Pourahmadi and Soofi (2000) and  
Pourahmadi (2001, Section 8.3).
However, those in (\ref{eq:ent}) have not appeared before.
For $n>0$, slightly more general calculations leading to analogues
of (\ref{eq:ent}) and (\ref{eq:4.15taking}) can be used to show
that the inverse autocorrelation function of $\{ X_t \}$ at lag
$u$ 
is the negative of the partial
correlation between $X_0$ and $X_u$ after elimination of the
effects of $X_t$, $t \neq 0, u$, as shown in Kanto (1984) for
processes with strictly positive spectral density functions. 
\end{ex}

\begin{ex}[The Finite Yaglom Problem]
There are many situations where the cardinality of $M$ is two or
more; see Pourahmadi {\it et al.}\ (2007), Box and Tiao (1975),
Brubacher and Wilson (1976), Damsleth (1980), Abraham (1981).
In the literature of time series analysis, there are several ad
hoc methods for interpolating the missing values.
For example, Brubacher and Wilson (1976) minimize
\[
\sum\nolimits^n_{-m} \varepsilon^2_j 
= \sum\nolimits^n_{-m} \left( \sum\nolimits^j_{k=- \infty}
a_{j-k} X_k \right)^2
\]
with respect to the unknown $X_j$, $j\in M\cup\{0\}$, and then
study the solution of the normal equations as
$m,n\rightarrow\infty$.
Budinsky (1989) has shown that this approach under some conditions
gives the same result as the more rigorous approach of Yaglom
(1963).
In applying Lemma \ref{L:snd} to this problem, we first note that,
due to the large cardinality of $M$, handling (\ref{eq:4.6aaa})
and (\ref{eq:4.5aaa}) via (\ref{eq:hhh}) does not lead to 
simple explicit formulas as in (\ref{eq:ent}) and
(\ref{eq:4.15taking}). 
Nevertheless, the limits of the expressions in (\ref{eq:4.6aaa})
and (\ref{eq:4.5aaa}) as first $m \rightarrow \infty$, and then as
$n \rightarrow\infty$ (assuming $f^{-1}\in L^1$) have simple forms
in terms of the AR parameters:
\begin{equation}
\left\{
\begin{aligned}
&X_0-\hat{X}_0(S)
=\sum\nolimits_{i\in M\cup\{0\}} \alpha^\prime_i
\left(
\sum\nolimits^\infty_{k=i}\bar{a}_{k-i}\,\varepsilon_k \right),\\
&\sigma^2(S)
=\mbox{the $(0,0)$-entry of the inverse of }
\left(\sum\nolimits^\infty_{k=i\vee j}
\bar a_{k-i}\,{a}_{k-j}\right)_{i,j\in M\cup\{0\}}.
\end{aligned}
\right.
\label{eq:4.16aaa}
\end{equation}

Now, using (\ref{eq:follows}) and writing the entries of the above
matrix, in terms of the Fourier coefficients of $f^{-1}$, it
follows that (\ref{eq:4.16aaa}) reduces to the results in  
Yaglom (1963); see also Salehi (1979). 

\end{ex}


\subsection{The Infinite Past and the Wold Decomposition}
\label{sec:InfPast}

A more direct method of solving prediction problems for  
$S=S_n\setminus M$ is to reduce them to a different class of
finite prediction problems than those in Section 3.2.
This is done by using the Wold decomposition of a purely
nondeterministic stationary process.

As in Section 3.1, write $N=\{-m,\ldots,n\}$ and
$N=K\cup\{0\}\cup M$ (disjoint), so that
$S=S_n\setminus M=\{\ldots,-m-2,-m-1\}\cup K$ (disjoint). 
For $j \geq -m$, let $\hat{X}_j$ be the linear least-squares
predictor of $X_j$  based on the infinite past $\{X_k;k<-m\}$. 
Then, by (\ref{eq:MA}), 
\[
X_j-\hat X_j
=\sum\nolimits_{k=-m}^jb_{j-k}\varepsilon_k,
\qquad
j\geq -m,
\]
which are orthogonal to $\overline{\mathrm{sp}}\{X_j; j<-m\}$, and
it follows that
\[
\overline{\mathrm{sp}}\{X_j; j\in S\}
=\mathrm{sp}\{X_j-\hat X_j; j\in K\}
\oplus\overline{\mathrm{sp}}\{X_j; j<-m\}.
\]
This equality plays the key role in finding the predictor of $X_0$
and its prediction  error variance, based on $\{X_j;j\in S\}$. 
In fact, by using it, we only have to solve the problem of
predicting $X_0-\hat X_0$ based on $\{X_j-\hat X_j;j\in K\}$. 
More precisely, we consider $X':=(X_j-\hat X_j)_{j\in N}$ which
has the covariance matrix $G=(g_{i,j})_{i,j\in N}$ with
\[
g_{i,j}:=\sum\nolimits_{k=-m}^{i\wedge j}b_{i-k}\bar b_{j-k}
\]
(see Pourahmadi (2001, p.\ 273)). 
Then, writing $X_0=\hat X_0+(X_0-\hat X_0)$, we get
\begin{equation}
\begin{split}
&\hat X_0(S)
=\hat X_0+\sum\nolimits_{k\in K}\alpha_k(X_k-\hat X_k),\\
&\sigma^2(S)=
\left\Vert(X_0-\hat X_0)
-\sum\nolimits_{k\in K}\alpha_k(X_k-\hat X_k)\right\Vert^2,
\end{split}
\label{eq:foo}
\end{equation}
where $\sum_{k\in K}\alpha_k(X_k-\hat X_k)$ is the predictor of
$X_0-\hat X_0$ based on $\{X_k-\hat X_k;k\in K\}$, and the
predictor coefficients $\alpha_k$ and prediction error variance
$\sigma^2(S)$ are obtained from the normal equations
(\ref{eq:normal}) with $\gamma_{i,j}$ replaced by $g_{i,j}$;
in particular, by (\ref{eq:sigma-minus2}),
\[
\sigma^2(S)
=\left[\mbox{the $(0,0)$-entry of the inverse of } 
\left(\sum\nolimits_{k=-m}^{i\wedge j}b_{i-k}\bar b_{j-k}\right)
_{i,j\in K\cup\{0\}}\right]^{-1}.
\]

We can also apply Lemma \ref{L:snd} to the above finite 
prediction problem for $X'$. 
In so doing, the following representations for 
the $(i,j)$-entry $g^{i,j}$ of $G^{-1}$ 
and the $j$-th entry $Y_j$ of the dual $Y$ of $X'$ are available:
\begin{equation}
g^{i,j}
=\sum\nolimits_{k=i\vee j}^n \bar{a}_{k-i}a_{k-j},
\qquad 
Y_j=\sum\nolimits_{k=j}^{n}\bar{a}_{k-j}\varepsilon_{k}.
\label{eq:4.17bbb}
\end{equation}
In fact, these are obtained by using (\ref{eq:ba}) and 
Proposition \ref{prop:3.2aaa} (3) 
or by letting $m\to\infty$ in (\ref{eq:4.4aaa}). 
The explicit representations in (\ref{eq:4.17bbb}) are also 
important in finding series representations for predictors and
interpolators discussed in the next two subsections.


\subsection{Series Representation of the Predictors} \label{sec:pred}

The Wold decomposition (\ref{eq:MA}) is often used to express
predictors and prediction errors in terms of the innovation
process $\{\varepsilon_t\}$. 
This strategy works well for achieving the goal (P$_2$) in
Section 1, but since the innovation $\varepsilon_t$ is not
directly observable the resulting predictor formulas are not
suitable for computation.
To get around this difficulty, one must express the innovations or
the predictors in terms of the past observations.
In this section, we obtain series representations for the infinite
past predictors in terms of the observed values. 
A novelty of our approach is its reliance on the
representation of the prediction error in terms of the dual $Y$
in (\ref{eq:4.17bbb}), hence the solution of the problem (P$_1$)
for $S=S_n\setminus M$ is more direct and simpler than the
procedures of Bondon (2002, Theorem 3.1) and Nikfar (2006).

Assuming that $\{X_j\}_{j\in\Z}$ has the mean-convergent AR
representation (\ref{eq:AR}),  it follows from (\ref{eq:foo}) with
$S=\{\ldots,-m-2,-m-1\}\cup K$ that
\[
\hat X_0(S)
=\sum\nolimits_{k\in K}\alpha_kX_k
+\sum\nolimits_{j=1}^\infty
\left(f_{j,m}-\sum\nolimits_{k\in K}\alpha_kf_{j,m+k}\right)
X_{-m-j},
\]
where $f_{j,k}:=-\sum_{i=0}^kb_{k-i}a_{j+i}$ is the coefficient of
the $(k+1)$-step ahead predictor based on the infinite past 
$S_0=\{\ldots,-2,-1\}$, i.e., 
$\hat X_k(S_0)=\sum_{j=1}^\infty f_{j,k}X_{-j}$
for $k=0,1,\dots$. 
On the other hand, from the finite duality principle or, more
precisely, (\ref{eq:fff}) with (\ref{eq:4.17bbb}), we have
\[
\hat X_0(S)
=X_0-\sum\nolimits_{i\in M\cup\{0\}}\alpha'_i
\left(\sum\nolimits_{k=i}^n\bar a_{k-i}\varepsilon_k\right).
\]
From this, replacing $\varepsilon_k$ from (\ref{eq:AR}) and after
some algebra, we get the following alternative series
representation for the predictor of $X_0$ based on the incomplete
past:
\begin{equation}
\hat X_0(S)
=-\sum\nolimits_{j\in S}
\left(\sum\nolimits_{i\in M\cup\{0\}}\alpha_i^\prime
\sum\nolimits_{k=i\vee j}^n\bar a_{k-i} a_{k-j}\right)X_j.
\label{eq:bondon}
\end{equation}
We note that the prediction error here has the representation
\begin{equation}
X_0-\hat X_0(S)
=\sum\nolimits_{i\in M\cup\{0\}}\alpha_i^\prime
\left(\sum\nolimits_{k=i}^n\bar a_{k-i}\varepsilon_k\right)
\label{eq:error_2}
\end{equation}
in terms of the dual $Y$ in (\ref{eq:4.17bbb}).
Furthermore, the sequence 
$\{\sum_{k=j}^{n}\bar a_{k-j}\varepsilon_k\}_{j=-\infty}^{n}$
spans $\overline{\mathrm{sp}}\{X_j; j\le n\}$, the infinite past
up to $n$ of the process $\{X_t\}$. 
The formulas (\ref{eq:bondon}) and (\ref{eq:error_2}) were
obtained initially by Bondon (2002, Theorem 3.2) without using the
notion of duality.


\subsection{Series Representation of the Interpolators}\label{sec:4.5aaa}

Series representation for the interpolator of $X_0$ based on the
observed values from the index set
$S=S_\infty\setminus M=\Z \backslash (M\cup\{0\})$
was obtained by Salehi (1979).
Here we obtain such representation using the idea of the dual
process.
Assuming $f^{-1}\in L^1$ or $\sum_{j=0}^\infty\abs{a_j}^2<\infty$,
the process 
\[
\xi_j:=\sum\nolimits_{k=j}^\infty\bar a_{k-j}\varepsilon_k,
\qquad
j\in\Z,
\]
is well-defined in the sense of mean-square convergence. 
From (\ref{eq:MA}), (\ref{eq:ba}), and the above results, 
we have the following: 
\begin{itemize}
\item[{\rm (i)}]
$(X_i,\xi_j)=\delta_{ij}$ for $i,j\in\Z$.
\item[{\rm (ii)}]
$\xi_j=\{X_j-\hat X_j(\Z\setminus\{j\})\}/
\Vert X_j-\hat X_j(\Z\setminus\{j\})\Vert^2$ for $j\in\Z$.
\item[{\rm (iii)}]
$\{\xi_j;j\in\Z\}$ spans the space
$\overline{\mathrm{sp}}\{X_j;j\in\Z\}$.
\item[{\rm (iv)}]
$\{ \xi_j; j \in \Z \}$ is a stationary process with the
autocovariance function
\[
\gamma^j
:=\frac1{2\pi}\int_{-\pi}^\pi
e^{-ij\lambda}f(\lambda)^{-1}d\lambda,
\qquad j\in\Z,
\]
i.e., 
$(\xi_i,\xi_j)=\gamma^{i-j}
=\sum_{k=i\vee j}^\infty\bar{a}_{k-i}a_{k-j}$ for $i,j\in\Z$.
\end{itemize}
The process $\{\xi_j\}$ has already appeared in prediction theory
and time series analysis, and is called the
{\it standardized two-sided innovation} (Masani (1960)) or the
{\it inverse process} (Cleveland (1972)) of $\{X_t\}_{t\in\Z}$.

Now, for solving the interpolation problem with
$S=\Z\setminus(M\cup\{0\})$, we need to show that
$\{\xi_j;j\in M\cup\{0\}\}$ spans the orthogonal complement of
$\overline{\mathrm{sp}}\{X_j; j\in S\}$ in
$\overline{\mathrm{sp}}\{X_j;j\in\Z\}$.
Then, it turns out that there is unique
$(\alpha_j^\prime)_{j\in M\cup\{0\}}$ satisfying
\[
X_0-\hat X_0(S)
=\sum\nolimits_{i\in M\cup\{0\}}\alpha_i^\prime\xi_i
=\sum\nolimits_{i\in M\cup\{0\}}\alpha_i^\prime
\left(\sum\nolimits_{k=i}^\infty\bar a_{k-i}\varepsilon_k\right)
\]
(see (\ref{eq:4.5aaa}) and (\ref{eq:error_2})),
and that $\sigma^2(S)=\alpha_0^\prime$.
Since
\[
(X_0,\xi_j)
-\sum\nolimits_{i\in M\cup\{0\}}\alpha_i^\prime(\xi_i,\xi_j)
=(\hat X_0(S),\xi_j)=0,
\qquad
j\in M\cup\{0\},
\]
we can compute $(\alpha_i^\prime)_{i\in M\cup\{0\}}$ by solving
the following system of linear equations:
\[
\sum\nolimits_{i\in M\cup\{0\}}\alpha_i^\prime \gamma^{i-j}
=\delta_{j0},
\qquad
j\in M\cup\{0\}.
\]
As for the predictor, if $\sum_{j=-\infty}^\infty \gamma^jX_{-j}$
is mean-convergent, then $(\xi_j)_{j\in\Z}$ admits the
representation
\[
\xi_i=\sum\nolimits_{j=-\infty}^\infty\gamma^{i-j}X_j,
\qquad
i\in\Z,
\]
and we have
\[
\hat X_0(S)
=-\sum\nolimits_{j\in S}
\left(\sum\nolimits_{i\in M\cup\{0\}}\alpha_i\gamma^{i-j}\right)
X_j,
\]
which is the two-sided version of the formula (\ref{eq:bondon}).


\section{Discussion and Future Work} \label{sec:discu}

We have reviewed and unified some important results from
prediction theory of stationary processes using a
finite-dimensional duality principle whose proof is based on
elementary ideas from the linear algebra.
Our time-domain, geometric and finite-dimensional approach brings
considerable clarity and simplicity to this area of prediction
theory as compared to the classical spectral-domain approach based
on analytic function theory and duality in the
infinite-dimensional spaces.
Since our duality lemma is not confined to stationary processes or
Toeplitz matrices, it has the potential of being useful in
solving similar prediction problems for nonstationary processes,
particularly those with low displacement ranks
(Kailath and Sayed (1995)).
However, the present form of the lemma does not seem to be useful
for prediction problems of infinite-variance or $L^p$-processes 
(Cambanis and Soltani (1984), Cheng {\it et al.}\ (1998)).


\end{document}